# Fast Decentralized Optimization over Networks

Meng Ma, *Student Member, IEEE*, Athanasios N. Nikolakopoulos, *Member, IEEE*,
and Georgios B. Giannakis, *Fellow, IEEE*

*Abstract*—The present work introduces the *hybrid consensus alternating direction method of multipliers* (H-CADMM), a novel framework for optimization over networks which unifies existing distributed optimization approaches, including the centralized and the decentralized consensus ADMM. H-CADMM provides a flexible tool that leverages the underlying graph topology in order to achieve a desirable *sweet-spot* between node-to-node communication overhead and rate of convergence – thereby alleviating known limitations of both C-CADMM and D-CADMM. A rigorous analysis of the novel method establishes linear convergence rate, and also guides the choice of parameters to optimize this rate. The novel hybrid update rules of H-CADMM lend themselves to "*in-network acceleration*" that is shown to effect considerable – and essentially "free-of-charge" – performance boost over the fully decentralized ADMM. Comprehensive numerical tests validate the analysis and showcase the potential of the method in tackling efficiently, widely useful learning tasks.

## I. INTRODUCTION

Recent advances in machine learning, signal processing and data mining, have led to important problems that can be formulated as distributed optimization over networks. Such problems entail parallel processing of data acquired by interconnected nodes and arise frequently in several applications, including data fusion and processing using sensor networks [45], [34], [24], [38]; vehicle coordination [37], [36]; power state estimation [23]; clustering [13]; classification [14]; regression [28]; filtering [33]; and demodulation [46], [1], to name a few. Among the candidate solvers for such problems, the *alternating direction method of multipliers* (ADMM) [3], [5] stands out as an efficient and easily implementable algorithm of choice that has attracted much interest in recent years [11], [18], [20], [7], thanks to its simplicity, fast convergence, and easily decomposable structure.

Many distributed optimization problems can be formulated in a consensus form, and solved efficiently by ADMM [5], [16]. The solver involves two basic steps: (i) a *communication step* for exchanging information with a central processing unit, the so-called *fusion center* (FC); and, (ii) an *update step* for updating the local variables at each node. By alternating between the two, local iterates eventually converge to the global solution. This approach is referred to as *centralized consensus ADMM* (C-CADMM), and although it has been successfully applied in various settings, it may not always present the preferable solver. In large-scale systems for instance, the cost of connecting each node to the FC may become prohibitive as the overhead of communicating data to the FC may be overwhelming, and the related storage requirement could surpass the capacity of a single FC. Furthermore, having one dedicated FC can lead to a single point of failure. In addition, there might be privacy-related issues that restrict access to private data.

Decentralized optimization on the other hand, forgoes with the FC by exchanging information only among single-hop neighbors. As long as the network is connected, local iterates can consent to the globally optimal decision variable, thanks to the aforementioned information exchange. This method – referred to as *decentralized* consensus ADMM (D-CADMM) – has attracted considerable interest; see e.g., [16] for a review of applications in communications and networking. In practice however, the fact that nodes only communicate with their neighbors may limit its convergence rate (especially as the network size increases), rendering it impractical for several large-scale applications.

### A. Our Contributions

To address the aforementioned limitations, the present paper puts forth a novel decentralized framework, that we term hybrid consensus ADMM, which unifies and markedly broadens C-CADMM and D-CADMM. Our contributions are in five directions:

(i) H-CADMM features *hybrid* updates accommodating communications with both the FCs and single-hop neighbors, thus bridging centralized with fully decentralized updates. This makes H-CADMM appealing for large-scale networks with multiple local FCs – a situation none of the existing approached is designed to handle.

(ii) A novel formulation of D-CADMM without duplicate constraints (dual variables commonly adopted by decentralized learning [23], [16], [39]) emerges simply by specializing the hybrid constraints to coincide with those arising from the purely neighborhood-based formulation.

(iii) Linear convergence is established, along with a rate bound, and specializes to C- and D-CADMM. The parameter setting to achieve the optimal bound is also provided.

(iv) H-CADMM is flexible to deploy FCs as needed to maximize performance gains, thus striking a desirable trade-off between the number of FCs deployed and convergence gain sought.

(v) The capability of handling hybrid constraints not only deals with mixed updates, but also effects "in-network acceleration" in decentralized operation without incurring noticeable increase in the overall complexity.

This work was supported in part by NSF grants 1442686, 1500713, 1509040, and 171141. The authors are with the Digital Technology Center and the Department of Electrical and Computer Engineering, University of Minnesota, Minneapolis, MN 55455, USA. Emails: maxxx971, anikolak, georgios@umn.edu.

## B. Related Work

Distributed optimization over networks has attracted much attention since the seminal works in [3], [44], where gradient-based parallel algorithms were developed. Since then, several alternatives have been advocated, including subgradient methods [29], [30], [41], [43], stochastic subgradients [40], dual averaging [9], [42], and gossip algorithms [8].

The decomposability of CADMM makes it particularly well suited for distributed optimization. Along with its many variants, including centralized [3], decentralized [16], [39], weighted [24], [27] and Nesterov accelerated [18], CADMM has gained wide popularity.

ADMM was introduced in the 1970s [17], and its convergence analysis can be found in e.g. [15], [10]. Local linear convergence of ADMM for linear or quadratic programs is established in [4]; see also [20] where the cost is a sum of component costs. Global linear convergence of a more general form of ADMM is reported in [7]; and linear convergence for a generalized formulation of consensus ADMM using the so-called "communication matrix" in [27].

Though D-CADMM has been applied to various problems [24], [46], [28], [38], [1], [2], its convergence remained open until linear convergence of D-CADMM was established in [39], and in [25] for its weighted counterpart. A successive orthogonal projection approach for distributed learning over networked nodes is introduced in [35], where nodes cannot communicate, but each node can access only limited amounts of data and agreement is enforced across nodes sharing the same data. A distributed ADMM algorithm that deals with node clusters was proposed, for which linear convergence was also substantiated [21]. However, it relies on the gossip algorithm for communication within clusters when cluster head is absent, and the explicit rate bound is difficult to obtain; and, it admits no closed-form representation for general networks. The present contribution is the first principled attempt to tackling the hybrid consensus problem.

## C. Outline and Notation

The rest of the paper is organized as follows. Section II states the problem and outlines two existing solvers, namely C-CADMM and D-CADMM; Section III develops H-CADMM, and shows its connections to both C-CADMM and D-CADMM; Section IV establishes linear convergence of H-CADMM, and discusses parameter settings that can afford optimal performance; Section V introduces the notion of "in-network acceleration;" Section VI reports the results of numerical tests; and Section VII concludes this work.

**Notation.** Vectors are represented by bold lower case, and matrices by bold upper case letters, $\mathbf{I}_N$ denotes the identity matrix of size $N \times N$ and $\mathbf{1}(\mathbf{0})$ the all ones (zeros) vector of appropriate size.

## II. PRELIMINARIES

For a network of $N$ nodes, consensus optimization amounts to solving problems of the form

$$\min_{\mathbf{x}} \quad \sum_{i=1}^{N} f_i(\mathbf{x}) \tag{1}$$

where $f_i(\cdot)$ is the $i$-th cost – only available to node $i$; and $\mathbf{x} \in \mathbb{R}^l$ is the common decision variable.

A common approach to solving such problems is to create a local copy of the global decision variable for each node, and impose equality constraints among all local copies; that is,

$$\min_{\{\mathbf{x}_i\}} \quad \sum_{i=1}^{N} f_i(\mathbf{x}_i) \tag{2}$$
$$\text{s. to} \quad \mathbf{x}_1 = \mathbf{x}_2 = \ldots = \mathbf{x}_N$$

where $\{\mathbf{x}_i\}$ are the local copies, and equality is enforced to ensure equivalence of (1) with (2). As a result, the global decision variable in (1), is successfully decoupled to facilitate distributed processing.

Each node optimizes locally its component of cost, and the equality constraints are effected by exchanging information among nodes, subject to restrictions. Indeed, in the *centralized* case nodes communicate with a single FC, while in the fully *decentralized* case, nodes can only communicate with their immediate neighbors.

We model communication constraints in the decentralized setting as an undirected graph $\mathcal{G} := (\mathcal{V}, \mathcal{E})$, with each vertex $\{v_i\}$ corresponding to one node, and the presence of edge $(v_i, v_j) \in \mathcal{E}$ denoting that nodes $i$ and $j$ can communicate. With $N$(resp. $M$) denoting the number of nodes (edges), we will label nodes (edges) using the set $\{1, 2, \ldots, N\}$ (resp. $\{1, 2, \ldots, M\}$). We will further define the neighborhood set of node $i$ as $\mathcal{N}_i := \{j | (v_i, v_j) \in \mathcal{E}\}$.

The following assumptions will be adopted about the graph and the local cost functions.

**Assumption 1** (Connectivity). Graph $\mathcal{G} := (\mathcal{V}, \mathcal{E})$ is connected.

**Assumption 2** (Strong convexity). Local cost $f_i$ is $\sigma_i$-strongly convex; that is, for any $\mathbf{x}, \mathbf{y} \in \mathbb{R}^l$,

$$f_i(\mathbf{y}) \geq f_i(\mathbf{x}) + \nabla^\top f_i(\mathbf{x})(\mathbf{y} - \mathbf{x}) + \frac{\sigma_i}{2}\|\mathbf{y} - \mathbf{x}\|_2^2.$$

**Assumption 3** (Lipschitz continuous gradient). Local $f_i$ is differentiable, and has Lipschitz continuous gradient; that is, for any $\mathbf{x}, \mathbf{y} \in \mathbb{R}^l$,

$$\|\nabla f_i(\mathbf{x}) - \nabla f_i(\mathbf{y})\|_2 \leq L_i\|\mathbf{x} - \mathbf{y}\|_2.$$

For brevity, we will henceforth focus on $l = 1$, but Appendix A outlines the generalization to $l \geq 2$.

## A. Centralized consensus ADMM

With a centralized global (G)FC, consensus is guaranteed when each node forces its local decision variable to equal that of the GFC. In iterative algorithms, this is accomplished through the update of each local decision variable based on information exchanged with the GFC. As a result, (1) can be formulated as

$$\min_{\{x_i\}} \quad \sum_{i=1}^{N} f_i(x_i) \tag{3}$$
$$\text{s. to} \quad x_i = z$$

where $z$ represents the GFC's decision variable (state).

ADMM solver by (3) by (i) forming the augmented Lagrangian; and (ii) performing Gauss-Seidel updates of primal and dual variables. Attaching Lagrange multipliers $\{\lambda_i\}_{i=1}^N$ to the equality constraints, and augmenting the Lagrangian with the penalty parameter $\rho$, we arrive at

$$L(\mathbf{x}, \mathbf{z}, \boldsymbol{\lambda}) = F(\mathbf{x}) + \boldsymbol{\lambda}^\top (\mathbf{x} - z\mathbf{1}) + \frac{\rho}{2}\|\mathbf{x} - z\mathbf{1}\|_2^2$$

where $\mathbf{x} := [x_1, x_2, \ldots, x_N]^\top$, $\boldsymbol{\lambda} := [\lambda_1, \lambda_2, \ldots, \lambda_N]^\top$, $F(\mathbf{x}) := \sum_{i=1}^N f_i(x_i)$, and $z = \bar{x} = N^{-1}\sum_{i=1}^N x_i$. Per entry (node) $i$, the ADMM updates are (see e.g. [5])

$$x_i^{k+1} = (\nabla f_i + I)^{-1}(\rho \bar{x}^k - \lambda_i^k) \tag{4a}$$

$$\lambda_i^{k+1} = \lambda_i^k + \rho(x_i^{k+1} - \bar{x}^k) \tag{4b}$$

where $z$ has been eliminated, and the inverse in (4a) is a shorthand for $x_i^{k+1}$ being the solution of $\nabla f_i(x_i^{k+1}) + x_i^{k+1} = (\rho \bar{x}^k - \lambda_i^k)$. Specialized to (4) C-CADMM boils down to the following three-step updates.

C1. Node $i$ solves (4a), and $x_i^{k+1}$ to the GFC;

C2. The GFC updates its global decision variable by averaging local copies, and broadcasts the updated value $z^{k+1}$ back to all nodes; and

C3. Node $i$ updates its Lagrange multiplier as in (4b).

### B. Decentralized consensus ADMM

In the decentralized setting, no GFC is present and nodes can only communicate with their one-hop neighbors. If the underlying graph $\mathcal{G}$ is connected, consensus constraints effect agreement across nodes.

Consider an auxiliary variable $\{z_{ij}\}$ per edge $(v_i, v_j) \in \mathcal{E}$, and re-write (2) in as

$$\min_{\{x_i\}} \quad \sum_{i=1}^N f_i(x_i) \tag{5}$$
$$\text{s. to} \quad x_i = z_{ij}, \; x_j = z_{ji}, \; (v_i, v_j) \in \mathcal{E}.$$

For undirected graphs, we have $z_{ij} = z_{ji}$. With $M$ edges in total, (5) includes $4M$ equality constraints, that can be written in matrix-vector form, leading to the compact expression

$$\min_{\{x_i\}} \quad F(\mathbf{x}) \tag{6}$$
$$\text{s. to} \quad \mathbf{A}\mathbf{x} + \mathbf{B}\mathbf{z} = \mathbf{0}$$

where $\mathbf{z}$ is the vector concatenating all $\{z_{ij}\}$ in arbitrary order, $\mathbf{A} := [\mathbf{A}_1^\top, \mathbf{A}_2^\top]^\top$ with $\mathbf{A}_1, \mathbf{A}_2 \in \mathbb{R}^{N \times 2M}$ defined such that if the $q$-th element of $\mathbf{z}$ is $z_{ij}$, then $(\mathbf{A}_1)_{qi} = 1$, $(\mathbf{A}_2)_{qj} = 1$, and all other elements are zeros; while $\mathbf{B} := [-\mathbf{I}_{2M}^\top, -\mathbf{I}_{2M}^\top]^\top$.

Formulation (6) is amenable to ADMM. To this end, one starts with the augmented Lagrangian

$$L(\mathbf{x}, \mathbf{z}, \boldsymbol{\lambda}) = F(\mathbf{x}) + \boldsymbol{\lambda}^\top (\mathbf{A}\mathbf{x} + \mathbf{B}\mathbf{z}) + \frac{\rho}{2}\|\mathbf{A}\mathbf{x} + \mathbf{B}\mathbf{z}\|_2^2$$

where Lagrange multiplier vector $\boldsymbol{\lambda} := [\boldsymbol{\beta}^\top, \boldsymbol{\gamma}^\top]^\top$ is split in sub-vectors $\boldsymbol{\beta}, \boldsymbol{\gamma} \in \mathbb{R}^{2M}$, initialized with $\boldsymbol{\beta}^0 = -\boldsymbol{\gamma}^0$.

After simple manipulations one obtains the simplified ADMM updates (see [39] for details):

$$x_i^{k+1} = (\nabla f_i + \rho |\mathcal{N}_i| I)^{-1}\left[\frac{\rho}{2}\sum_{j \in \mathcal{N}_i}(x_i^k + x_j^k) - y_i^k\right] \tag{7a}$$

$$y_i^{k+1} = y_i^k + \frac{\rho}{2}\sum_{j \in \mathcal{N}_i}\left(x_i^{k+1} - x_j^{k+1}\right) \tag{7b}$$

where $\mathbf{y} := (\mathbf{A}_1 - \mathbf{A}_2)^\top \boldsymbol{\beta}$, and the inverse in (7a) is a shorthand for $x_i^{k+1}$ being the solution of $\nabla f_i(x_i^{k+1}) + \rho |\mathcal{N}| x_i^{k+1} = (\rho/2)\sum_{j \in \mathcal{N}_i}(x_i^{k+1} + x_j^{k+1})$.

The fact that the per-node updates in (7) involve only single-hop neighbors justifies the term decentralized consensus ADMM (D-CADMM).

In a nutshell, D-CADMM works as follows:

D1. Each node sends its local variable to all its single-hop neighbors;

D2. Upon receiving information from all its neighbors, node $i$ updates its local variable as in (7a);

D3. Node $x_i$, node $i$ updates its dual variable $y_i$ as in (7b).

## III. Hybrid consensus ADMM

Rather than a single GFC that is connected to all nodes, here we consider optimization over networks with *multiple LFCs*. Such a setup can arise in large-scale networks, where bandwidth, power, and computational limits or even security concerns may discourage deployment of a single GFC. These considerations prompt the deployment of multiple LFCs each of which communicates with a limited number of nodes. No prior ADMM-based solver can deal with this setup as none is capable of handling *hybrid constraints* that are present when nodes exchange information not only with LFCs but also with their single-hop neighbors. This section, introduces H-CADMM that is particularly designed to handle this situation.

### A. Problem formulation

In contrast to the simple graph $\mathcal{G}$ used in the fully decentralized setting, we will employ *hypergrahs* to cope with hybrid constraints. A hypergraph is a tuple $\mathcal{H} := (\mathcal{V}, \mathcal{E})$, where $\mathcal{V}$ is the vertex set and $\mathcal{E} = \{\mathcal{E}_i\}_{i=1}^M$ denotes the collection of hyperedges. Each $\mathcal{E}_i$ comprises a set of vertices, $\mathcal{E}_i \subset \mathcal{V}$ with cardinality $|\mathcal{E}_i| \geq 2$, $\forall i$. A vertex $v_i$ and an edge $\mathcal{E}_j$ are said to be incident if $v_i \in \mathcal{E}_j$. Hypergraphs are particularly suitable for modeling hybrid constraints because each LFC can be modeled as one hyperedge consisting of all its connected nodes.

With $N$ nodes, $M$ hyperedges, and their corresponding orderings, we can associate each edge variable $z_j$ with hyperedge $j$. Then the hybrid constraints can be readily reparameterized as $x_i = z_j$, $\forall i : v_i \in \mathcal{E}_j$. Consider now vectors $\mathbf{x} \in \mathbb{R}^N$, $\mathbf{z} \in \mathbb{R}^M$ collecting all local $\{x_i, z_j\}$s, and matrices $\mathbf{A} \in \mathbb{R}^{T \times N}$, $\mathbf{B} \in \mathbb{R}^{T \times M}$ constructed to have nonzero entries $A_{ti} = 1$, $B_{tj} = 1$ corresponding to $t$-th constraint $x_i - z_j = 0$. For $T$ equality constraints, the hybrid form of (1) can thuse be written compactly as

$$\min_{x_i} \quad \sum_{i=1}^N f_i(x_i) \tag{8}$$
$$\text{s.to} \quad \mathbf{A}\mathbf{x} - \mathbf{B}\mathbf{z} = \mathbf{0}.$$

Let now $\mathbf{C} \in \mathbb{R}^{N \times M}$ denote the incidence matrix of the hypergraph, formed with entries $C_{ij} = 1$ if node $i$ and edge $j$ are incident, and $C_{ij} = 0$ otherwise; $d_i$ the degree of node $i$ (the number of incident edges of node $i$); $e_j$ the degree of hyperedge $j$ (the number of incident nodes of hyperedge $j$); diagonal matrix $\mathbf{D} \in \mathbb{R}^{N \times N}$ the node degree matrix (formed with $d_i$ as its $i$-th diagonal element); and likewise $\mathbf{E} \in \mathbb{R}^{M \times M}$ the *edge degree* matrix (formed with $e_j = |\mathcal{E}_j|$ as its $j$-th diagonal element). With these notational conventions, we prove in the Appendix the following.

**Lemma 1.** *Matrices $\mathbf{A}$ and $\mathbf{B}$ in (8) satisfy*

$$\mathbf{A}^{\top} \mathbf{A} = \mathbf{D} \tag{9a}$$

$$\mathbf{B}^{\top} \mathbf{B} = \mathbf{E} \tag{9b}$$

$$\mathbf{A}^{\top} \mathbf{B} = \mathbf{C}. \tag{9c}$$

### B. Algorithm

The augmented Lagrangian for (8) is

$$L(\mathbf{x}, \mathbf{z}, \boldsymbol{\lambda}) = \sum_{i=1}^{N} f_i(x_i) + \boldsymbol{\lambda}^{\top}(\mathbf{A}\mathbf{x} - \mathbf{B}\mathbf{z}) + \frac{\rho}{2}\|\mathbf{A}\mathbf{x} - \mathbf{B}\mathbf{z}\|_2^2 \tag{10}$$

where $\boldsymbol{\lambda} \in \mathbb{R}^T$ collects all the Lagrange multipliers, and $\rho$ is a hyper-parameter controlling the effect of the quadratic regularizer. ADMM updates can be obtained by cyclically solving for $\mathbf{x}$, $\mathbf{z}$ and $\boldsymbol{\lambda}$ the equations

$$\nabla f(\mathbf{x}^{k+1}) + \mathbf{A}^{\top} \boldsymbol{\lambda}^k + \rho \mathbf{A}^{\top}(\mathbf{A}\mathbf{x}^{k+1} - \mathbf{B}\mathbf{z}^k) = \mathbf{0} \tag{11a}$$

$$\mathbf{B}^{\top} \boldsymbol{\lambda}^k + \rho \mathbf{B}^{\top}(\mathbf{A}\mathbf{x}^{k+1} - \mathbf{B}\mathbf{z}^{k+1}) = \mathbf{0} \tag{11b}$$

$$\boldsymbol{\lambda}^{k+1} = \boldsymbol{\lambda}^k + \rho(\mathbf{A}\mathbf{x}^{k+1} - \mathbf{B}\mathbf{z}^{k+1}). \tag{11c}$$

Equations (11) can be simplified by left-multiplying (11c) by $\mathbf{B}^{\top}$ and adding it to (11b) to obtain

$$\mathbf{B}^{\top} \boldsymbol{\lambda}^{k+1} = \mathbf{0}. \tag{12}$$

If $\boldsymbol{\lambda}$ is initialized such that $\mathbf{B}^{\top} \boldsymbol{\lambda}^0 = \mathbf{0}$, then $\mathbf{B}^{\top} \boldsymbol{\lambda}^k = \mathbf{0}$ for all $k \geq 0$. Eliminating $\mathbf{B}^{\top} \boldsymbol{\lambda}^k$ from (11b), one can solve for $\mathbf{z}$, and arrive at the closed form

$$\mathbf{z}^{k+1} = \mathbf{E}^{-1} \mathbf{C}^{\top} \mathbf{x}^{k+1}. \tag{13}$$

Similarly, by left-multiplying (11c) by $\mathbf{A}^{\top}$ and letting $\mathbf{y}^k := \mathbf{A}^{\top} \boldsymbol{\lambda}^k$, one finds

$$\mathbf{y}^{k+1} - \mathbf{y}^k = \rho(\mathbf{D}\mathbf{x}^k - \mathbf{C}\mathbf{z}^{k+1}). \tag{14}$$

Then simply plugging (13) into (11a) yields

$$\mathbf{x}^{k+1} = (\nabla f + \rho \mathbf{D} I)^{-1}(\rho \mathbf{C}\mathbf{z}^k - \mathbf{y}^k). \tag{15}$$

Recursions (13)–(15) summarize our H-CADMM, and their per-node forms are listed in Algorithm 1.

*Remark* 1. Two interesting observations are in order:

(i) Regardless of the number of attached nodes, each hyperedge serves as an LFC; and

(ii) Each LFC performs local averaging. Indeed, the entry-wise update of (13) shows that each hyperedge satisfies $z_i = (1/E_{ii}) \sum_{j \in \mathcal{N}_i} x_j$. Hence, all hyperedges are treated equally in the sense that they are updated by the average value of all incident nodes.

---

**Algorithm 1:** Hybrid Consensus ADMM

**Input:** $\rho$, $\mathbf{x}^0$, $\mathbf{z}^0$, $\mathbf{y}^0 = \mathbf{0}$
**while** *stopping criterion not satisfied* **do**
  **for** $i = 1, \ldots, N$ **do**
    node $i$ updates $x_i^{k+1}$ by solving
    $\nabla f_i(x_i^{k+1}) + \rho|\mathcal{N}_i|x_i^{k+1} = \rho \sum_{j \in \mathcal{N}_i} z_j^k - y_i^k$
    send $x_i^{k+1}$ to all incident FCs and neighbors
  **for** $j = 1, \ldots, M$ **do**
    FC $j$ updates $z_j^{k+1} = \frac{1}{E_{jj}} \sum_{i \in \mathcal{N}_j} x_i^{k+1}$
    send $z_j^{k+1}$ to all incident nodes
  **for** $i = 1, \ldots, N$ **do**
    node $i$ updates
    $y_i^{k+1} = y_i^k + \rho\left(D_{ii}x_i^{k+1} - \sum_{j \in \mathcal{N}_i} z_j^{k+1}\right)$

---

### C. Key relations

Here we unveil a relationship satisfied by the iterates $\{\mathbf{x}^k\}$ generated by Algorithm 1, which not only provides a different view of H-CADMM, but also serves as the starting point for establishing the convergence results in Section IV. This relation shows that $\mathbf{x}^{k+1}$ depends solely on the gradient of the local cost function, as well as the past $\{\mathbf{x}^k, \mathbf{x}^{k-1}, \ldots, \mathbf{x}^0\}$ that is also not dependent on the variables $\mathbf{z}^k$ and $\mathbf{y}^k$.

**Lemma 2.** *The sequence $\{\mathbf{x}^k\}$ generated by Algorithm 1 satisfies*

$$\mathbf{x}^{k+1} = -\frac{1}{\rho}\mathbf{D}^{-1}\nabla F(\mathbf{x}^{k+1}) + \mathbf{D}^{-1}\mathbf{C}\mathbf{E}^{-1}\mathbf{C}^{\top}\mathbf{x}^k$$
$$- \mathbf{D}^{-1}(\mathbf{D} - \mathbf{C}\mathbf{E}^{-1}\mathbf{C}^{\top})\sum_{t=0}^{k} \mathbf{x}^t. \tag{16}$$

*Proof.* Substituting (13) into (14), we obtain

$$\mathbf{y}^{k+1} - \mathbf{y}^k = \rho(\mathbf{D} - \mathbf{C}\mathbf{E}^{-1}\mathbf{C}^{\top})\mathbf{x}^k \tag{17}$$

which upon initializing with $\mathbf{y}^0 = \mathbf{0}$, leads to

$$\mathbf{y}^k = \rho(\mathbf{D} - \mathbf{C}\mathbf{E}^{-1}\mathbf{C}^{\top})\sum_{t=0}^{k} \mathbf{x}^t. \tag{18}$$

Plugging (13) and (18) into (15), yields

$$\nabla F(\mathbf{x}^{k+1}) + \rho \mathbf{D}\mathbf{x}^{k+1} = \rho \mathbf{C}\mathbf{E}^{-1}\mathbf{C}^{\top}\mathbf{x}^k$$
$$- \rho(\mathbf{D} - \mathbf{C}\mathbf{E}^{-1}\mathbf{C}^{\top})\sum_{t=0}^{k} \mathbf{x}^t \tag{19}$$

from which we can readily solve for $\mathbf{x}^{k+1}$. $\qquad\square$

Lemma 2 shows that $x_i^{k+1}$ is determined by its past $\{x_i^t\}_{t=0}^k$, and the local gradient, namely $\nabla f_i(x_i^{k+1})$. This suggests a new update scheme, where each node maintains not only its current $x_i^k$ but also $\sum_{t=0}^{k} x_i^t$.

Among the things worth stressing in Lemma 2 is the appearance of matrices $\mathbf{C}\mathbf{E}^{-1}\mathbf{C}^{\top}$ and $\mathbf{D} - \mathbf{C}\mathbf{E}^{-1}\mathbf{C}^{\top}$. Since both play key roles in studying the evolution of (19), it is important to understand their properties and impact on the performance of the algorithm.

**Lemma 3.** *Matrices* $\mathbf{CE}^{-1}\mathbf{C}^\top$ *and* $\mathbf{D} - \mathbf{CE}^{-1}\mathbf{C}^\top$ *are positive semidefinite (PSD), and satisfy*

$$(\mathbf{D} - \mathbf{CE}^{-1}\mathbf{C}^\top)\mathbf{1} = \mathbf{0}. \tag{20}$$

*Proof.* See Appendix. $\qquad\square$

### D. H-CADMM links to C-CADMM and D-CADMM

Modeling hybrid communication constraints as a hypergraph not only affords the flexibility to accommodate multiple LFCs, but also provides a unified view of consensus-based ADMM. Indeed, it is not difficult to show that by specializing the hypergraph, our proposed approach subsumes both centralized and decentralized consensus ADMM.

**Proposition 1.** *H-CADMM reduces to C-CADMM when there is only one hyperedge capturing all nodes.*

*Proof.* When there is a single hyperedge comprising all network nodes, we have $\mathbf{A} = \mathbf{I}_N$, and $\mathbf{B} = \mathbf{1}$. Using (9), we thus obtain

$$\mathbf{D} = \mathbf{A}^\top\mathbf{A} = \mathbf{I}_N \tag{21a}$$

$$\mathbf{E} = \mathbf{B}^\top\mathbf{B} = N \tag{21b}$$

$$\mathbf{C} = \mathbf{A}^\top\mathbf{B} = \mathbf{1}. \tag{21c}$$

Then, the update (13) at the GFC reduces to

$$z^{k+1} = \mathbf{E}^{-1}\mathbf{C}^\top\mathbf{x}^k = \frac{1}{N}\sum_{i=1}^N x_i^k = \bar{x}^k.$$

Similarly, (15) and (14) specialize to

$$\mathbf{x}^{k+1} = (\nabla f + \rho I)^{-1}(\rho\bar{x}^k\mathbf{1} - \boldsymbol{\lambda}^k) \tag{22a}$$

$$\boldsymbol{\lambda}^{k+1} = \boldsymbol{\lambda}^k + \rho(\mathbf{x}^{k+1} - \bar{x}^{k+1}\mathbf{1}) \tag{22b}$$

where we have used that $\mathbf{y}^k = \mathbf{A}^\top\boldsymbol{\lambda}^k = \boldsymbol{\lambda}^k$.

Comparing (22) with (4), it is not difficult to see that (4a) is just the entry-wise form of (22a); and likewise for (22b) and (4b). Therefore C-CADMM can be viewed as a special case of H-CADMM with one hyperedge connecting all nodes. $\qquad\square$

**Proposition 2.** *H-CADMM reduces to D-CADMM when every edge is a hyperedge.*

*Proof.* When each simple edge is modeled as a hyperedge, the resulting hyperedges will end up having degree 2, that is, $\mathbf{E} = 2\mathbf{I}$. Thus, (13) becomes

$$\mathbf{z}^{k+1} = \frac{1}{2}\mathbf{C}^\top\mathbf{x}^{k+1}. \tag{23}$$

Using (23), and eliminating $\mathbf{z}$ from (15) and (14) yields

$$\mathbf{x}^{k+1} = (\nabla f + \rho\mathbf{D}I)^{-1}(\frac{\rho}{2}\mathbf{CC}^\top\mathbf{x}^k - \mathbf{y}^k) \tag{24a}$$

$$\mathbf{y}^{k+1} = \mathbf{y}^k + \rho(\mathbf{D} - \frac{1}{2}\mathbf{CC}^\top)\mathbf{x}^{k+1}. \tag{24b}$$

To relate $\sum_{j\in\mathcal{N}}(x_i + x_j)$ in (7b) to (24b), notice that $D_{ii} = |\mathcal{N}_i|$ and $(\mathbf{CC}^\top\mathbf{x})_i = \sum_{j\in\mathcal{N}_i}(x_i + x_j)$. Hence, it is straightforward to see that (7a) and (7b) are just the entry-wise versions of (24a) and (24b). Therefore, D-CADMM is also a special case of H-CADMM with each simple edge viewed as a hyperedge. $\qquad\square$

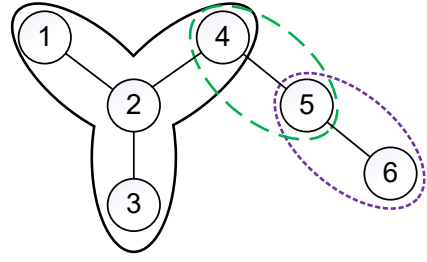

Fig. 1. The simple graph of Example 1. Both the underlying graph (edges in black lines) and the hypergraph (hyperedges as ellipsoids) are shown for comparison.

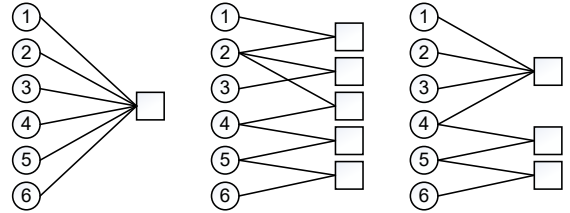

Fig. 2. Communication graphs of C-CADMM, D-CADMM and H-CADMM for Example 1. Circles represent nodes, while squares represent hyperedges. Solid lines between nodes and hyperedges indicate nodes belonging to hyperedges.

*Remark* 2. In past works of fully decentralized consensus ADMM [28], [46], [16], [39], one edge is often associated with two variables $z_{ij}$, $z_{ji}$ in order to decouple the equality constraint $x_i = x_j$, and express it as $x_i = z_{ij}$, $x_j = z_{ji}$. Although eventually the duplicate variables are shown equal (and therefore discarded), their sheer presence leads to duplicate Lagrange multipliers, which can complicate the algorithm. Proposition 2 suggests a novel derivation of D-CADMM with only one variable attached to each edge, a property that can simplify the whole process considerably.

**Example 1.** Consider the simple graph depicted in Fig. 1 with 6 nodes and 5 edges. All nodes work collectively to minimize some separable cost. Solid lines denote the undirected graph connectivity, while ellipsoids represent hyperedges in the modeling hypergraph of H-CADMM. Let us consider C-CADMM, D-CADMM and H-CADMM solvers on this setting. To gain insight, we examine closely the update rules at node 4 that we list in Table I.

As C-CADMM relies on a GFC connecting all nodes, every node receives updates from the GFC. This is demonstrated in the first row of Table I.

D-CADMM on the other hand, allows communication only along edges (solid line). Thus, each node can only receive updates from its single-hop neighbors. Specifically, node 4 can only receive information from nodes 2 and 5, as can be seen in the second row of Table I.

TABLE I
COMPARISON OF ADMM UPDATE RULES AT NODE 4 USING C-CADMM, D-CADMM AND H-CADMM TO SOLVE THE PROBLEM IN EXAMPLE 1

| Method | Update Rules |
|---|---|
| **C-CADMM** | $x_4^{k+1} = \arg\min_{x_4} \left[ f_4(x_4) + \lambda_4^k(x_4 - \bar{x}^k) + \frac{\rho}{2}\|x_4 - \bar{x}^k\|_2^2 \right]$ |
| | $\lambda_4^{k+1} = \lambda_4^k + \rho(x_4^{k+1} - \bar{x}^{k+1})$ |
| **D-CADMM** | $x_4^{k+1} = (\nabla f_4 + 2\rho I)^{-1}\left( \rho x_4^k + \frac{\rho}{2}(x_2^k + x_5^k) - y_4^k \right)$ |
| | $y_4^{k+1} = y_4^k + \frac{\rho}{2}\left( 2x_4^{k+1} - x_2^{k+1} - x_5^{k+1} \right)$ |
| **H-CADMM** | $x_4^{k+1} = \arg\min_{x_4} \left[ f_4(x_4) + y_4^k x_4 + \frac{\rho}{2}\left( \|x_4 - \frac{1}{4}\sum_{i=1}^{4} x_i^k\|^2 + \|x_4 - \frac{1}{2}\sum_{i=1}^{2} x_i^k\|^2 \right) \right]$ |
| | $y_4^{k+1} = y_4^k + \rho\left( 2x_4^{k+1} - \frac{1}{4}\sum_{i=1}^{4} x_i^{k+1} - \frac{1}{2}\sum_{i=4}^{5} x_i^{k+1} \right)$ |

H-CADMM lies somewhere in between. It allows deployment of multiple LFCs, each of which is connected to a subset of nodes. The hypergraph model is shown in Figure 1 with hyperedges marked by ellipsoids. The union of 3 solid ellipsoids corresponds to the LFC, while dash and dotted ones represent the edges between nodes 4, 5, and 6. Speaking of node 4, the information needed to update its local variables comes from its neighbors, the LFC represented by $\sum_{i=1}^{4} x_i^{k+1}/4$; and node 5 represented by $(x_4^{k+1} + x_5^{k+1})/2$. That is exactly what we see in the third row of Table I.

*Remark 3.* For all three consensus ADMM algorithms, the *fusion centers* – both global and local – act as averaging operators that compute, store and broadcast the mean values of the local estimates from all connected nodes.

*Remark 4.* The dual update per node acts as an accumulator forming the sum of residuals between the node and its connected *fusion centers*. With $\mathbf{L} := 2(\mathbf{D} - \mathbf{C}\mathbf{E}^{-1}\mathbf{C}^\top)$ denoting the Laplacian of the hypergraph, the dual update per node boils down to

$$\mathbf{y}^{k+1} = \mathbf{y}^k + \frac{\rho}{2}\mathbf{L}\mathbf{x}^{k+1}, \quad \forall k \geq 1.$$

If dual variables are initialized such that $\mathbf{y}^0 = \mathbf{0}$, then

$$\mathbf{y}^k = \frac{\rho}{2}\mathbf{L}\sum_{t=1}^{k}\mathbf{x}^t.$$

This observation holds for all three algorithms.

*Remark 5.* Figure 2 exemplifies that H-CADMM "lies" somewhere between C-CADMM and D-CADMM. Clearly, consensus is attainable if and only if the communication graph is connected.

## IV. CONVERGENCE RATE ANALYSIS

In this section, we analyze the convergence behavior of the novel H-CADMM algorithm. In particular, our main theorem establishes linear convergence and provides a bound on the rate of convergence, which depends on properties of both the objective function, as well as the underlying graph topology.

Apart from the assumptions made in Section II, here we also need an additional one:

**Assumption 4.** There exists at least one saddle point $(\mathbf{x}^\star, \mathbf{z}^\star, \mathbf{y}^\star)$ of Algorithm 1 that satisfies the KKT conditions:

$$\nabla f(\mathbf{x}^\star) + \mathbf{y}^\star = \mathbf{0} \tag{25a}$$
$$\mathbf{z}^\star - \mathbf{E}^{-1}\mathbf{B}^\top \mathbf{A}\mathbf{x}^\star = \mathbf{0} \tag{25b}$$
$$(\mathbf{I} - \mathbf{B}\mathbf{E}^{-1}\mathbf{B}^\top)\mathbf{A}\mathbf{x}^\star = \mathbf{0}. \tag{25c}$$

This assumption is required for the development of Algorithm 1 as well as for the analysis of its convergence rate. If (as4) does not hold, either the original problem is unsolvable, or, it entails unbounded subproblems, or, a diverging sequence of $\boldsymbol{\lambda}^k$ [7].

Assumptions 1–4 guarantee the existence of at least one optimal solution. In fact, we can further prove that any saddle point is actually the unique solution of the KKT conditions (25), and hence of Algorithm 1.

**Lemma 4.** *If $\boldsymbol{\lambda}$ is initialized so that $\mathbf{B}^\top\boldsymbol{\lambda}^0 = \mathbf{0}$, and (as1)–(as4) hold, then $(\mathbf{x}^\star, \mathbf{z}^\star, \mathbf{y}^\star)$ is the* unique *optimal solution of* (25).

*Proof.* See Appendix. □

### A. Linear rate of convergence

Alternating direction methods, including ADMM, have been thoroughly investigated in [7]. Similar to D-CADMM [39], conditions for establishing linear convergence rate in [7] are not necessarily satisfied by the H-CADMM setup[1]. Therefore, we cannot establish linear convergence rate simply by reformulating it as a special case of existing ADMM approaches.

One way to overcome this obstacle is to adopt a technique similar to [39], as we did in [26], to obtain a relatively loose bound on convergence, in the sense that it could not capture significant accelerations observed in practice by varying the topology of the underlying graph. For strongly convex costs, a tighter bound has been reported recently [27]. However, H-CADMM is not amenable to the analysis in [27] since the linear constraint coefficients $\mathbf{A}$ and $\mathbf{B}$ cannot be recovered from the communication matrix.

---

[1]This should be expected since H-CADMM reduces to D-CADMM upon modeling each simple edge as an hyperedge.

We establish convergence by measuring the progress in terms of the G-norm i.e. the semi-norm defined by $\|\mathbf{x}\|_G := \mathbf{x}^\top \mathbf{G} \mathbf{x}$, where $\mathbf{G}$ is the PSD matrix,

$$\mathbf{G} = \begin{bmatrix} \mathbf{I} & \mathbf{0} \\ \mathbf{0} & \mathbf{CE}^{-1}\mathbf{C}^\top \end{bmatrix}. \tag{26}$$

The G-norm is properly defined since both $\mathbf{CE}^{-1}\mathbf{C}^\top$ and $\mathbf{D} - \mathbf{CE}^{-1}\mathbf{C}^\top$ are PSD (see Lemma 3). Consider now the square root, $\mathbf{Q} := (\mathbf{D} - \mathbf{CE}^{-1}\mathbf{C}^\top)^{1/2}$, and define two auxiliary sequences

$$\mathbf{r}^k := \sum_{t=0}^{k} \mathbf{Q}\mathbf{x}^t, \quad \mathbf{q}^k := \begin{bmatrix} \mathbf{r}^k \\ \mathbf{x}^k \end{bmatrix}. \tag{27}$$

These two sequences play an important role in establishing linear convergence of the proposed algorithm. Before we establish such a convergence result we first need to bound the gradient of $F(\cdot)$.

**Lemma 5.** *If (as1)–(as4) hold, then for any $k \geq 0$, we have*

$$\mathbf{CE}^{-1}\mathbf{C}^\top(\mathbf{x}^{k+1} - \mathbf{x}^k) = -\mathbf{Q}(\mathbf{r}^{k+1} - \mathbf{r}^\star) - \frac{1}{\rho}\left(\nabla F(\mathbf{x}^{k+1}) - \nabla F(\mathbf{x}^\star)\right). \tag{28}$$

*Proof.* See Appendix. $\square$

Let $\Lambda := \lambda_N(\mathbf{CE}^{-1}\mathbf{C}^\top)$ denote the largest eigenvalue of $\mathbf{CE}^{-1}\mathbf{C}^\top$ and $\lambda := \lambda_2(\mathbf{D} - \mathbf{CE}^{-1}\mathbf{C}^\top)$ the second smallest eigenvalue of $\mathbf{D} - \mathbf{CE}^{-1}\mathbf{C}^\top$.

**Theorem 1.** *Under (as1)–(as4) for any $\rho > 0$, $\beta \in (0, 1)$, and $k > 0$, H-CADMM iterates in (8) satisfy*

$$\|\mathbf{x}^k - \mathbf{x}^\star\|_G^2 \leq \left(\frac{1}{1+\delta}\right)^k \|\mathbf{q}^0 - \mathbf{q}^\star\|_G^2 \tag{29}$$

*where $\mathbf{G}$ and $\mathbf{q}$ as in (26) and (27), and $\delta$ satisfies*

$$\delta \leq \min\left\{ \frac{2\beta\sigma}{\rho(\Lambda + \frac{2\Lambda}{\lambda})}, \frac{(1-\beta)\rho\lambda}{L} \right\}. \tag{30}$$

*Proof.* See Appendix. $\square$

Theorem 1 asserts that $\mathbf{x}^k$ converges linearly to the optimal solution $\mathbf{x}^\star$ at a rate bounded by $1/(1+\delta)$. Larger $\delta$ implies faster convergence.

Note that while Theorem 1 is proved for $l = 1$, it can be generalized to $l \geq 2$ (see (40) in Appendix A).

### B. Fine-tuning the parameters

Theorem 1 characterizes the convergence of iterates generated by H-CADMM. Parameter $\delta$ is determined by the local costs, the underlying communication graph topology, and the scalar $\rho$. By tuning these parameters, one can maximize the convergence bound, to speed up convergence in practice too. With local costs and the graph fixed, one can maximize $\delta$ by tuning $\rho$. When possible to choose the number and locations of LFCs, we can effectively alter the topology of the communication graph – hence modify $\Lambda$ and $\lambda$ – to improve convergence.

**Theorem 2.** *Under assumptions 1–4 the optimal convergence rate bound*

$$\delta \leq \frac{1}{\sqrt{\frac{L}{\sigma}\frac{\Lambda}{\lambda}(1 + 2\frac{\Lambda}{\lambda})}} \tag{31}$$

*is achieved by setting*

$$\rho = \sqrt{\frac{2\sigma L}{\Lambda\lambda(1 + \frac{\Lambda}{\lambda})}}. \tag{32}$$

*Proof.* The optimal $\beta^\star \in (0, 1)$ maximizing $\delta$ in Theorem 1 is

$$\beta^\star = \frac{\rho^2 \Lambda\lambda(1 + 2\frac{\Lambda}{\lambda})}{2\sigma L + \rho^2 \Lambda\lambda(1 + 2\frac{\Lambda}{\lambda})} \tag{33}$$

and is obtained by equating the two terms in (30)

$$\frac{2\beta\sigma}{\rho(\Lambda + 2\frac{\Lambda}{\lambda})} = \frac{(1-\beta)\rho\lambda}{L}. \tag{34}$$

Substituting (23) into (34), one arrives at

$$\delta = \frac{2\sigma\rho\lambda}{2\sigma L + \rho^2 \Lambda\lambda(1 + 2\frac{\Lambda}{\lambda})}. \tag{35}$$

Maximizing $\delta$ by varying $\rho$ eventually leads to (31), with the optimal

$$\rho^\star = \sqrt{\frac{2\sigma L}{\Lambda\lambda(1 + 2\frac{\Lambda}{\lambda})}}. \tag{36}$$

$\square$

Upon defining the cost condition number as

$$\kappa_F := \max_i \frac{L_i}{\sigma_i}$$

and the graph condition number as

$$\kappa_G := \frac{\lambda_N(\mathbf{CE}^{-1}\mathbf{C}^\top)}{\lambda_2(\mathbf{D} - \mathbf{CE}^{-1}\mathbf{C}^\top)} \tag{37}$$

the optimal convergence rate can be bounded as

$$\delta \leq \frac{1}{\sqrt{\kappa_F \kappa_G(1 + 2\kappa_G)}}. \tag{38}$$

Clearly, the bound in (38) is a decreasing function of $\kappa_F$ and $\kappa_G$. Therefore, decreasing both will drive the bound larger, possibly resulting in a faster rate of convergence. On the one hand, smaller cost condition number makes the cost easier to optimize; on the other hand, smaller graph condition number implies improved connectivity. Indeed, when the communication hypergraph is just a simple graph – a case for which H-CADMM reduces to D-CADMM with $\mathbf{D} - \mathbf{CE}^{-1}\mathbf{C}^\top = \mathbf{L}/2$, then $\lambda$ is the smallest nonzero eigenvalue of the Laplacian, which is related to bottlenecks in the underlying graph [6].

*Remark 6.* Theorem 2 shows that the number of iterations it takes to achieve an $\epsilon$-accurate solution is $\mathcal{O}(\sqrt{\kappa_F}\log(\frac{1}{\epsilon}))$. The dependence on $1/\sqrt{\kappa_F}$ improves over [39], which had established a dependence of $1/\kappa_F^2$.

## V. Graph-aware acceleration

Distributed optimization over networks using a central GFCs is not feasible for several reasons including communication constraints and privacy concerns. At the other end of the spectrum, decentralized schemes relying on single-hop communications may suffer from slow convergence, especially when the network has a large diameter or bottlenecks. H-CADMM fills the gap by compromising between the two aforementioned extremes. By carefully deploying multiple LFCs, it becomes possible to achieve significant performance gains whilst abiding by cost and privacy constraints.

In certain cases, leveraging the topology of the LFCs deployed could bring sufficient gains. Instead of, or complementing gains from these *actual* LFCs, this section advocates that gains in H-CADMM convergence are possible through *virtual FCs* on judiciously selected nodes. We refer to the benefit brought by virtual FCs as *in-network acceleration* (see Figure 3 for a simple illustration); it will be confirmed by numerical tests, virtual FCs can afford a boost in performance essentially "for free"; simply by exploiting the actual network topology.

The merits of in-network acceleration through virtual FCs at a subset of selected nodes (*hosts*), can be recognized in the following four aspects.

- *Hardware*. Relying on virtual LFCs, in-network acceleration requires no modifications in the actual topology and hardware.
- *Interface*. The other nodes "see" exactly the same number of neighbors, so there is no change in the communication interface. However, the information exchanged is indeed different.
- *Computational complexity*. Except for the host nodes, the update rules for both primal and dual variables remain the same. Each host however, serves a dual role: as an FC, as well as an ordinary node. We know from Algorithm 1 that the computations performed per LFC involve averaging information from all connected nodes, which is simple compared to updating the local variables. Thus, the computational complexity remains of the same order, while the total computational cost decreases as less iterations are necessary to reach a target level of accuracy.
- *Communication cost*. Since there is no change of the communication interface, the communication cost remains invariant. Once again, the total communication cost can further drop, since H-CADMM enjoys faster convergence.

Given that the interface does not change, nor extra communication/computation cost is incurred, one can think of in-network acceleration as a sort of "free lunch" approach, with particularly attractive practical implications.

### A. Strategies for selecting the local FCs

A reasonable question to ask at this point is: *"How should one select the nodes to host the virtual FCs?"* Unfortunately, there is no simple answer. The question would have been easier if we could choose as many LFCs as necessary to achieve the maximum possible acceleration. In practice however, we

---

**Algorithm 2:** Greedy Selection of LFCs

**Input:** LFC budget $B$, graph $\mathcal{G} := (\mathcal{V}, \mathcal{E})$
$\mathcal{H} \leftarrow$ empty set
$c \leftarrow 0$
**while** $\mathcal{V}$ *not empty* **do**
  mark $v \in \mathcal{V}$ with largest node degree as a LFC
  add hyperedge $\{v\} \cup \mathcal{N}_v$ to $\mathcal{H}$
  delete $\{v\} \cup \mathcal{N}_v$ from $\mathcal{V}$
  delete $\{(v_1, v_2) | v_1, v_2 \in \{v\} \cup \mathcal{N}_v\}$ from $\mathcal{E}$
  $c \leftarrow c + 1$
  **if** $c \geq B$ **then**
    break
$\mathcal{H} \leftarrow \mathcal{H} \cup \mathcal{E}$
generate $\mathbf{C}$ from $\mathcal{H}$
**Output:** incidence matrix $\mathbf{C}$

---

do not always have the luxury to place as many virtual LFCs as we want, for reasons that include lack of control over some nodes, and difficulty to modify internal updating rules. And even if we could, picking the right nodes to host the LFCs under a general network architecture might not be straightforward.

A reasonable way forth would be to maximize the convergence rate bound, $\delta$, subject to a maximum number of LFCs, hoping that the optimal solution would yield the best rate of convergence in practice. However, this turns out to involve optimizing the ratio of eigenvalues, which is typically difficult to solve. For this reason, we will resort to heuristic methods.

Intuitively, one may choose the nodes with highest degree so as to maximize the effect of virtual LFCs. However, one should be careful when applying this simple approach to clustered graphs. For example, consider the graph consisting of two cliques (connected by a short path), comprising $n_1$ and $n_2$ nodes, respectively. Each node in the larger clique has higher degree than every node in the smaller one. As a result, always assigning the role of LFC to the largest degree nodes would disregard the nodes of the smaller clique (when our budget is less than $n_1$), while one could apparently take care of both cliques with as few as two LFCs. Taking this into account, we advocate a greedy LFC selection (Algorithm 2), which prohibits placing virtual LFCs within the neighborhood of other FCs.

*Remark 7.* For simplicity, Algorithm 2 relies on degree information to select LFCs. More elaborate strategies would involve richer structural properties of the underlying topology, to identify more promising nodes at the expense of possibly computationally heavier LFC selection. In general, LFC selection strategies offer the potential of substantially increasing the convergence rate over random assignment. Note however, that regardless of the choice of virtual LFCs, H-CADMM remains operational.

### B. On H-CADMM's "free lunch"

In-network acceleration has several advantages, allowing for faster convergence essentially "without paying any price." At first glance, this appears to be a "free lunch" type of benefit, and deservedly makes one skeptical. Actually, the benefit comes from leveraging information that is completely overlooked by fully decentralized methods. To see this, recall that in D-CADMM, each node communicates with only one neighbor each time, without accounting for the entire neighborhood. Instead, H-CADMM manages to exploit network topology by creating virtual LFCs that gather and share information with the whole neighborhood. It is this additional information that enables faster flow of data, and hence faster convergence. Therefore, it is not a "free lunch" for H-CADMM; but a lunch "not even tasted" by D-CADMM.

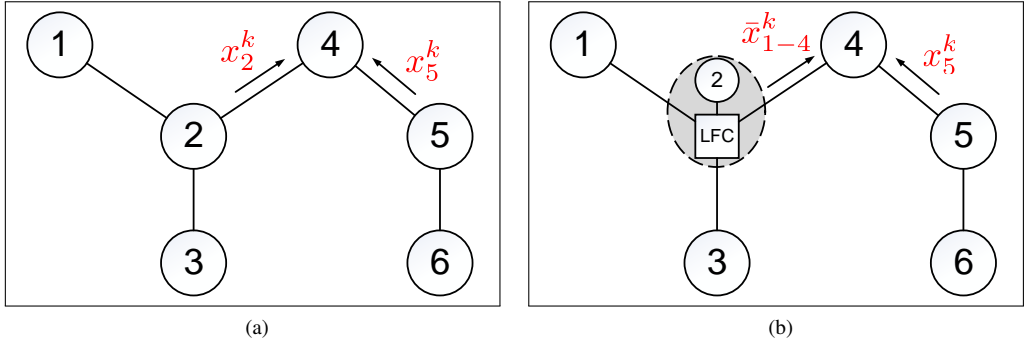

Fig. 3. A demonstration of in-network acceleration applied to the problem of Example 1. The shaded dashed circle in (b) is equivalent to node 2 in (a), except that a virtual FC (square LFC) is created logically, making it amenable to the application of H-CADMM. The interface to other nodes remains the same. The information exchanged through the edges changes. For example, the message sent from node 2 to 4 changes from $x_2^k$, which describes only the state of node 2, to $x_{1-4}^k = 1/N \sum_{n=1}^4 x_n^k$ which contains information about 4 nodes.

## VI. NUMERICAL TESTS

In this section, we test numerically the performance of H-CADMM, and also validate our analytical findings.

### A. Experimental settings

Throughout this section, we consider several interconnected nodes trying to estimate one value, $x_0$, based on local observations $o_i = x_0 + \epsilon_i$, where $\epsilon_i \sim \mathcal{N}(0, 0.1)$ is the measurement noise. This can be solved by minimizing the least-squares (LS) error $F(\mathbf{x}) = \frac{1}{2} \sum_{i=1}^N \|o_i - x_i\|_2^2$. Different from centralized LS, here the observation $o_i$ is available only to node $i$, and all nodes collaborate to obtain the final solution. In D-CADMM each node can only talk to its neighbors, while in H-CADMM nodes can potentially communicate with the LFC and their neighbors.

We test D-CADMM and H-CADMM solvers with various parameter settings. We assess convergence using the relative accuracy metric defined as $\|x^k - x^\star\|_2 / \|x^\star\|_2$, and report the number of iterations as well as the communication cost involved in reaching a target level of performance. The communication cost measures how many times local and global decision variables are transferred across the network. Originally, we set $\rho$ according to Theorem 2, but this choice did not work well our tests. For this reason, we tuned it manually to reach the best possible performance.

### B. Acceleration of dedicated FCs

In this test, we compare the performance of H-CADMM with dedicated FCs against that of D-CADMM. In particular, we choose only one dedicated LFC connected to 20% and 50% of the nodes drawn randomly from (i) a lollipop graph; (ii) a caveman graph; and (iii) two Erdos Renyi random graphs. All the graphs have $N = 50$ nodes. Specifically, in the lollipop graph, 50% of nodes comprise a clique and the rest form a line graph attached to this clique. The caveman graph consists of a cycle formed by 10 small cliques, each forming a complete graph of 5 nodes. The Erdos Renyi graphs are randomly generated with edge probability $r = 0.05$ and $r = 0.1$, respectively.

Figure 4 compares the performance of H-CADMM with one dedicated FC connected to 20% and 50% of the nodes against D-CADMM, in terms of number of iterations needed to achieve certain accuracy, as well as, communication cost. Lines with the same color markers denote the same graph; dashed lines correspond to D-CADMM, and solid lines to H-CADMM. From these two figures, one can draw several interesting observations:

- H-CADMM with mixed updates works well in practice. Solutions are obtained in fewer iterations, and at lower communication cost.
- The tests verify the linear convergence properties of both D-CADMM and H-CADMM, as can be seen in all four graphs.
- The performance gap between dashed lines and solid lines of the same marker confirms the acceleration ability of H-CADMM. The gap is larger for "badly-connected" graphs, such as the lollipop graph, and relatively small for "well-connected" graphs, such as the Erdos Renyi graphs. In fact, this observation holds even when comparing the Erdos Renyi graphs. Indeed, for the ER($r = 0.05$) graph which is not as well connected as the ER($r = 0.1$), the performance gap is smaller.
- Figure 4 suggests that the more nodes are connected to the FC, the larger the acceleration gains for H-CADMM, which is intuitively reasonable since extra connections pay off. In view of this connections-versus-acceleration trade-off, H-CADMM can reach desirable sweet spots

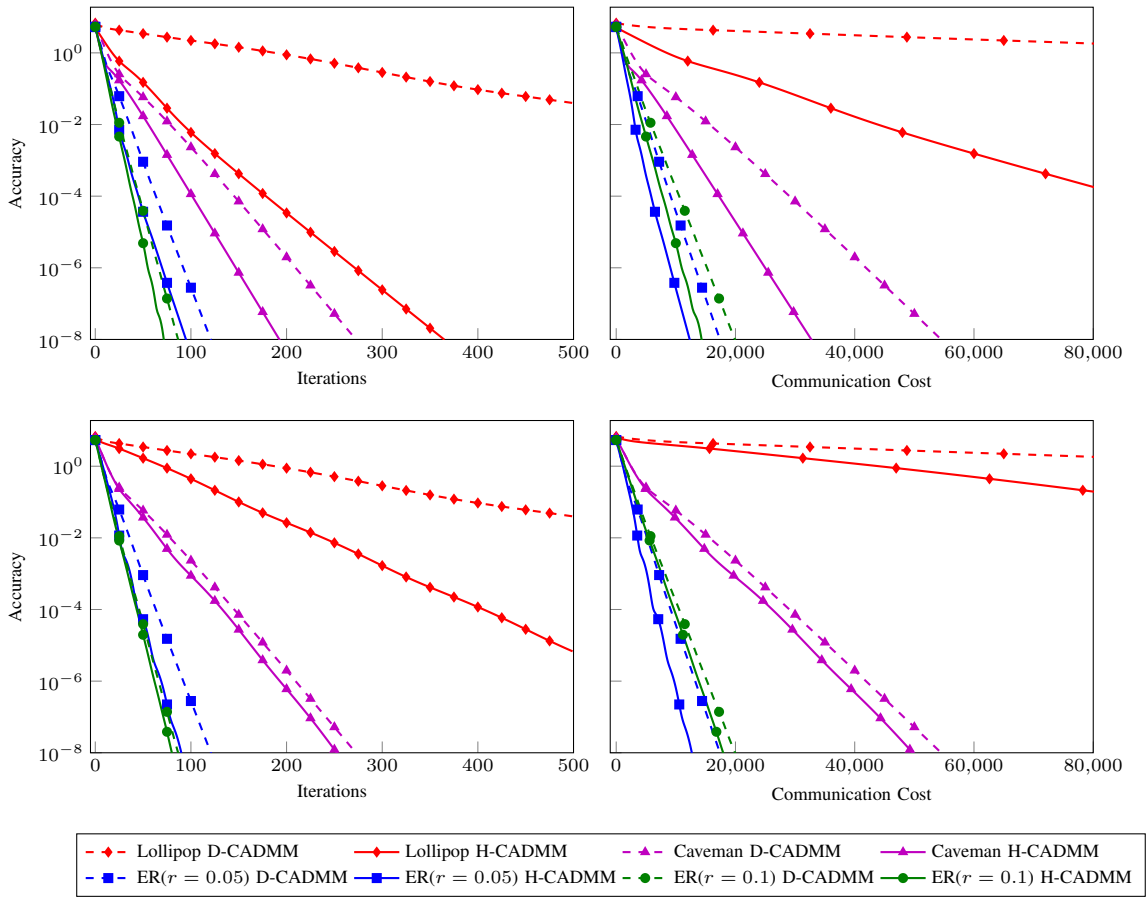

Fig. 4. Performance comparison of D-CADMM and H-CADMM in terms of iteration number as well as communication cost. H-CADMM is configured with one dedicated FC connecting 50% of the nodes (top row figures) and 20% of the nodes respectively (bottom row figures).

between performance gains and deployment cost.

### C. In-network acceleration

In this test, we demonstrate the performance gain effected by in-network acceleration. By creating virtual FCs among nodes, this technique does not require dedicated FCs and new links. As detailed in Section V, in-network accelerated H-CADMM exchanges information along existing edges, essentially leading to a communication cost that follows the same pattern with iteration complexity. Therefore, in this experiment, we report both metrics in one figure.

We first apply H-CADMM with in-network acceleration to several fixed-topology graphs, namely the line graph, the cycle graph, and the star graph, and we report the results in Figure 5. Then, we carry out the same tests on the lollipop graph, the caveman graph and the two Erdos-Renyi graphs with parameters $r = 0.05$ and $r = 0.10$, and we present the results in Figure 6. In each test, we select the hosts using Algorithm 2.

All the results illustrate that H-CADMM with in-network acceleration offers a significant boost in convergence rate over D-CADMM, especially for graphs with relatively large diameters (or graphs that are not well-connected), such as the line graph, the cycle graph or the lollipop graph. On the other hand, the performance gain is minimal for the star graph, whose diameter is 2 regardless of the number of nodes, as well as the Erdos Renyi random graph with high edge probability (see Figure 6). Note that these performance gains over D-CADMM are achieved without paying a substantial computational cost (just one averaging step), which speaks for the practical merits of H-CADMM.

### D. Trade-off between FCs and performance gain

Finally, we explore the trade-off between the number of LFCs and the corresponding convergence rate. We perform tests on several graphs with different properties, and using only H-CADMM with in-network acceleration. In this test, we measure performance by the number of iterations needed to achieve a target accuracy of $10^{-8}$, given a varying number of LFCs ranging from 1 to 25.

Figure 7 depicts the results. In general, as the number of LFCs increases, the number of iterations decreases. Besides this general trend, one can also make the following observations.

- For Erdos Renyi graphs with a relatively high edge probability, there is a small initial gain arising from the introduction of virtual LFCs, which diminishes fast as their number increases. This is not surprising since such graphs are "well connected," and therefore, adding more virtual LFCs does not help much. On the other hand, for "badly connected" graphs, such as the line graph or the lollipop graph, there is a significant convergence boost as the number of LFCs increases.

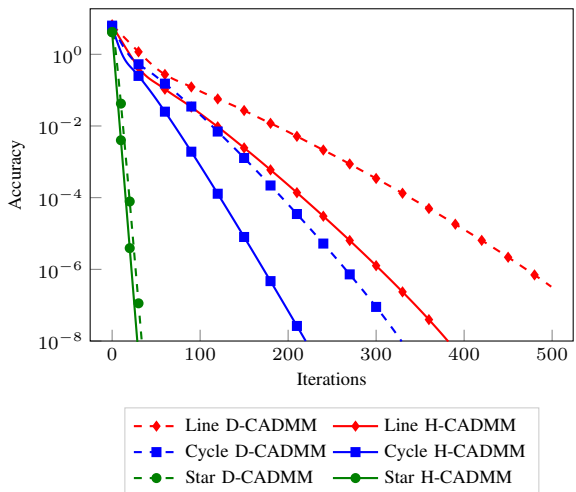

Fig. 5. Effects of in-network acceleration on the line graph, the cycle graph, and the star graph.

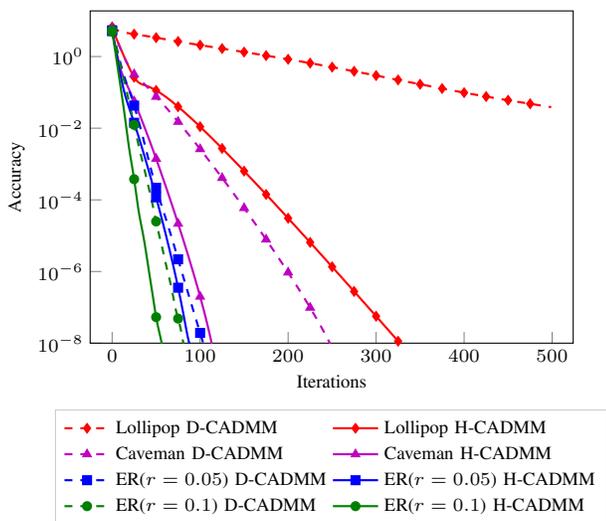

Fig. 6. Effects of in-network acceleration on the lollipop graph, the caveman graph, and two Erdos-Renyi random graphs.

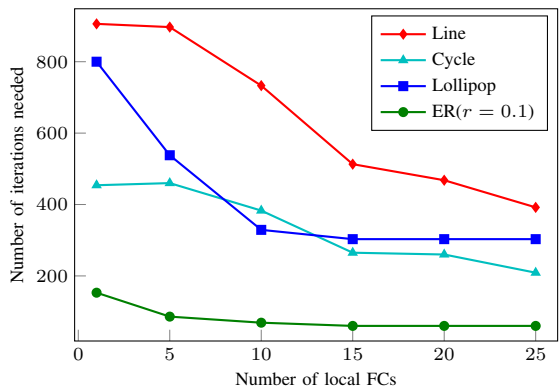

Fig. 7. The impact of adding LFCs in H-CADMM with in-network acceleration. The performance is measured in terms of number of iterations needed to achieve certain accuracy.

- For the line and the cycle graph, a significant change in convergence rate happens only after an initial threshold has been surpassed (in our case 5-6 LFCs). For the lollipop graph, there seems to exist a cut-off point above which adding more nodes does not lead to significant change in convergence rate.

## VII. Conclusion and Future Work

This paper introduces the novel H-CADMM algorithm that generalizes the centralized and the decentralized CADMM, while also accelerating D-CADMM with modified updates.

We establish linear convergence of H-CADMM, and we also conduct a comprehensive set of numerical tests that validate our analytical findings and demonstrate the effectiveness of the proposed approach in practice.

A very promising direction we are currently pursuing involves the development of techniques leveraging the intrinsic hierarchical organization that is commonly found in distributed system architectures (see e.g., [19], [32]) as well as real-world large-scale network topologies (see e.g. [12], [31]).

## Appendix A
## Algorithm 1 for $l > 2$

For $l \geq 2$, we have $\mathbf{x} \in \mathbb{R}^{Nl}$, $\mathbf{z} \in \mathbb{R}^{Ml}$, $\boldsymbol{\lambda} \in \mathbb{R}^{Tl}$, and $\mathbf{y} \in \mathbb{R}^{Nl}$. Let $\tilde{\mathbf{A}} \in \mathbb{R}^{Tl \times Nl}$ and $\tilde{\mathbf{B}} \in \mathbb{R}^{Tl \times Ml}$ denote the coefficient matrices for $l \geq 2$, obtained by replacing 1's of $\mathbf{A}$ and $\mathbf{B}$ with the identity matrix $\mathbf{I}_l$ and 0's with all-zero matrix $\mathbf{0}_{l \times l}$. Consequently, node degree, edge degree, and incidence matrices are $\tilde{\mathbf{D}} = \tilde{\mathbf{A}}^\top \tilde{\mathbf{A}} \in \mathbb{R}^{Nl \times Nl}$, $\tilde{\mathbf{E}} = \tilde{\mathbf{B}}^\top \tilde{\mathbf{B}} \in \mathbb{R}^{Ml \times Ml}$, and $\tilde{\mathbf{C}} = \tilde{\mathbf{A}}^\top \tilde{\mathbf{B}} \in \mathbb{R}^{Nl \times Ml}$. Meanwhile, H-CADMM boils down to

$$\mathbf{x}^{k+1} = (\nabla f + \rho \tilde{\mathbf{D}} I)^{-1} (c \tilde{\mathbf{C}} \mathbf{z}^k - \mathbf{y}^k) \tag{39a}$$

$$\mathbf{z}^{k+1} = \tilde{\mathbf{E}}^{-1} \tilde{\mathbf{C}}^\top \mathbf{x}^{k+1} \tag{39b}$$

$$\mathbf{y}^{k+1} = \mathbf{y}^k + \rho (\tilde{\mathbf{D}} \mathbf{x}^{k+1} - \tilde{\mathbf{C}} \mathbf{z}^{k+1}). \tag{39c}$$

Relative to (13)–(15), the computational complexity of (39) is clearly higher. To see this, consider the per-step complexity for $l = 1$. Given that $\mathbf{E}$ is a diagonal matrix, naive matrix vector multiplication incurs complexity $\mathcal{O}(MN)$, which can be reduced to $\mathcal{O}(ME_{\max})$ by exploiting the sparsity of $\mathbf{C}$, where $E_{\max}$ is the largest edge degree. When $l \geq 2$ however, per-step complexity grows to $\mathcal{O}(l^2 M E_{\max})$ which – being quadratic in $l$ – would make it difficult for the method to handle high-dimensional data. Thankfully, a compact form of (39) made possible by Proposition 3 reduces the complexity from quadratic to linear in $l$.

**Proposition 3.** *Let $\mathbf{X} \in \mathbb{R}^{N \times l}$ denote the matrix formed with $i$-th row $\mathbf{x}_i^\top$ and likewise for $\mathbf{Z} \in \mathbb{R}^{M \times l}$ and $\mathbf{Y} \in \mathbb{R}^{N \times l}$. Then,* (39) *is equivalent to*

$$\mathbf{X}^{k+1} = (\nabla f + \rho \mathbf{D} I)^{-1} (\rho \mathbf{C} \mathbf{Z}^k - \mathbf{Y}^k) \tag{40a}$$

$$\mathbf{Z}^{k+1} = \mathbf{E}^{-1} \mathbf{C}^\top \mathbf{X}^{k+1} \tag{40b}$$

$$\mathbf{Y}^{k+1} = \mathbf{Y}^k + \rho (\mathbf{D} \mathbf{X}^{k+1} - \mathbf{C} \mathbf{Z}^{k+1}). \tag{40c}$$

*Proof.* The block structure suggests the following compact representation using Kronecker products

$$\tilde{\mathbf{C}} = \mathbf{C} \otimes \mathbf{I}_d, \quad \tilde{\mathbf{D}} = \mathbf{D} \otimes \mathbf{I}_d, \quad \tilde{\mathbf{E}} = \mathbf{E} \otimes \mathbf{I}_d.$$

Exploiting properties of Kronecker products [22, §2.8], matrices $\mathbf{P}$, $\mathbf{Q}$, $\mathbf{R}$, and $\mathbf{S}$ with compatible dimensions, satisfy

$$(\mathbf{P} \otimes \mathbf{Q}^\top)\mathbf{x} = \text{vec}(\mathbf{PXQ})$$

where $\mathbf{x}$ is obtained by concatenating all rows of $\mathbf{X}$ that we denote as $\mathbf{x} = \text{vec}(\mathbf{X}^\top)$. Thus, by setting $\mathbf{P} = \mathbf{C}$ and $\mathbf{Q} = \mathbf{I}_d$, one arrives at

$$\tilde{\mathbf{C}}\mathbf{z} = (\mathbf{C} \otimes \mathbf{I}_d)\mathbf{z} = \text{vec}(\mathbf{CZI}_d) = \text{vec}(\mathbf{CZ})$$
$$\tilde{\mathbf{D}}\mathbf{x} = (\mathbf{D} \otimes \mathbf{I}_d)\mathbf{x} = \text{vec}(\mathbf{DXI}_d) = \text{vec}(\mathbf{DX})$$

from which one readily obtains (40a) and (40c). To see the equivalence of (39b) and (40b), we use another property of Kronecker products, namely

$$(\mathbf{P} \otimes \mathbf{Q})(\mathbf{R} \otimes \mathbf{S}) = (\mathbf{PR}) \otimes (\mathbf{QS}).$$

Since $\tilde{\mathbf{E}}$ is block diagonal, so is $\tilde{\mathbf{E}}^{-1} = \mathbf{E}^{-1} \otimes \mathbf{I}_d$. Therefore, (39b) is equivalent to

$$\tilde{\mathbf{E}}^{-1}\tilde{\mathbf{C}}^\top \mathbf{x} = (\mathbf{E}^{-1} \otimes \mathbf{I}_d)(\mathbf{C}^\top \otimes \mathbf{I}_d)\text{vec}(\mathbf{X})$$
$$= ((\mathbf{E}^{-1}\mathbf{C}^\top) \otimes \mathbf{I}_d)\text{vec}(X) = \text{vec}(\mathbf{E}^{-1}\mathbf{C}^\top \mathbf{X})$$

which concludes the proof. $\square$

Proposition 3 establishes that H-CADMM can run using much smaller matrices, effectively reducing its computational complexity. To see this, consider the per-step complexity of (40). The difference with (39) is dominated by $\mathbf{C}^\top \mathbf{X}$, which leads to complexity $\mathcal{O}(lMN)$. This can be further reduced to $\mathcal{O}(lME_{\max})$ by exploiting the sparsity of $\mathbf{C}$, thus improving the complexity of (39) by a factor of $l$.

## Appendix B
## Proof of Lemma 1

Let $\mathbf{a}_i$ denote the $i$-th column of $\mathbf{A}$, and $\mathbf{b}_j$ the $j$-th column of $\mathbf{B}$. By construction, the $i$-th column of $\mathbf{A}$ indicates in which constraint $x_i$ is present; hence, $\mathbf{a}_i^\top \mathbf{a}_i = d_i$ equals the degree of node $i$. Similarly, it follows that $\mathbf{b}_j^\top \mathbf{b}_j = e_j$, and

$$\mathbf{a}_i^\top \mathbf{a}_j = 0, \quad \forall i \neq j, \qquad \mathbf{b}_i^\top \mathbf{b}_j = 0, \quad \forall i \neq j$$

from which we obtain $\mathbf{A}^\top \mathbf{A} = \mathbf{D}$ and $\mathbf{B}^\top \mathbf{B} = \mathbf{E}$.

Consider now the dot product $\mathbf{a}_i^\top \mathbf{b}_j$. When the $t$-th constraint reads $x_i = z_j$, it holds by construction that $(\mathbf{a}_i)_t = 1$, and $(\mathbf{b}_j)_t = 1$. Therefore, if node $i$ and edge $j$ are incident, we have $\mathbf{a}_i^\top \mathbf{b}_j = 1$; otherwise, $\mathbf{a}_i^\top \mathbf{b}_j = 0$. Thus, the incidence matrix definition implies that $\mathbf{A}^\top \mathbf{B} = \mathbf{C}$.

## Appendix C
## Proof of Lemma 3

Let $e_{\max}$ ($e_{\min}$) denote the maximum (minimum) degree of all edges, and $c_{ij}$ the number of common edges between nodes $i$ and $j$. Clearly, we have $d_i = \sum_{j=1}^{N} c_{ij}$, with $c_{ii} = 0$; and since $\mathbf{E} \preceq e_{\max}\mathbf{I}$, we obtain

$$\mathbf{CE}^{-1}\mathbf{C}^\top \succeq \frac{1}{e_{\max}}\mathbf{CC}^\top \succeq \mathbf{0}$$

where the last inequality holds because $\mathbf{C}$ has linearly independent columns, and hence $\mathbf{CC}^\top$ is PSD. The latter implies that $\mathbf{CE}^{-1}\mathbf{C}^\top$ is PSD too.

The definition of $\mathbf{C}$ leads to

$$\mathbf{CC}^\top = \begin{bmatrix} d_1 & c_{12} & \dots & c_{1N} \\ c_{21} & d_2 & \dots & c_{2N} \\ \vdots & \vdots & \ddots & \vdots \\ c_{N1} & c_{N2} & \dots & d_N \end{bmatrix}.$$

And since $\mathbf{E} \succeq e_{\min}\mathbf{I}$, it holds that

$$\mathbf{D} - \mathbf{CE}^{-1}\mathbf{C}^\top \succeq \mathbf{D} - \frac{1}{e_{\min}}\mathbf{CC}^\top \succeq \mathbf{D} - \frac{1}{2}\mathbf{CC}^\top$$
$$= \frac{1}{2}\begin{bmatrix} d_1 & -c_{12} & \dots -c_{1N} \\ -c_{21} & d_2 & \dots -c_{2N} \\ \vdots & \vdots & \ddots & \\ -c_{N1} & -c_{N2} & \dots & d_N \end{bmatrix}$$

where the last matrix is a valid Laplacian ($d_i = \sum_{j=1}^{N} c_{ij}$), which in turn implies that it is PSD.

Finally, by definition $\mathbf{C1} = \mathbf{d}$, $\mathbf{C}^\top \mathbf{1} = \mathbf{e}$, where $\mathbf{d} := [d_1, \dots, d_N]^\top$, $\mathbf{e} := [e_1, \dots, e_N]^\top$; and hence we have

$$(\mathbf{D} - \mathbf{CE}^{-1}\mathbf{C}^\top)\mathbf{1} = \mathbf{d} - \mathbf{C1} = \mathbf{0}.$$

## Appendix D
## Proof of Lemma 4

Feasibility of (25c) guarantees the existence of at least one solution, and therefore there exists at least one optimal solution. The uniqueness of $\mathbf{x}^\star$ and $\mathbf{z}^\star$ follows from the strong convexity of $F(\cdot)$ and the dual feasibility (25c).

To see the uniqueness of $\mathbf{y}$, we first show that if $\tilde{\boldsymbol{\lambda}}$ is optimal, then $\boldsymbol{\lambda}^\star$, the projection of $\tilde{\boldsymbol{\lambda}}$ to the column space of $\mathbf{A}$, is also optimal. Using the orthogonality principle, we arrive at

$$\mathbf{A}^\top(\tilde{\boldsymbol{\lambda}} - \boldsymbol{\lambda}^\star) = \mathbf{0}$$

which implies that $\boldsymbol{\lambda}^\star$ also satisfies (25a). Thus, projection of any solution $\tilde{\boldsymbol{\lambda}}$ to the column space of $\mathbf{A}$ is also an optimal solution.

If there are two optimal solutions, $\boldsymbol{\lambda}_1 \neq \boldsymbol{\lambda}_2$ in the column space of $\mathbf{A}$, we have

$$\mathbf{A}^\top \boldsymbol{\lambda}_1 = \mathbf{A}^\top \boldsymbol{\lambda}_2. \tag{41}$$

Furthermore, there exist $\mathbf{x}_1 \neq \mathbf{x}_2$ such that $\boldsymbol{\lambda}_1 = \mathbf{Ax}_1$, and $\boldsymbol{\lambda}_2 = \mathbf{Ax}_2$. Subtracting $\mathbf{A}^\top \boldsymbol{\lambda}_2$ from $\mathbf{A}^\top \boldsymbol{\lambda}_1$ yields

$$\mathbf{A}^\top \boldsymbol{\lambda}_1 - \mathbf{A}^\top \boldsymbol{\lambda}_2 = \mathbf{A}^\top \mathbf{A}(\mathbf{x}_1 - \mathbf{x}_2) = \mathbf{0} \tag{42}$$

from which we find that $\mathbf{x}_1 = \mathbf{x}_2$, which is a contradiction. Therefore, $\boldsymbol{\lambda}_1 = \boldsymbol{\lambda}_2$, and thus $\mathbf{y} = \mathbf{A}^\top \boldsymbol{\lambda}$ is also unique.

## APPENDIX E
## PROOF OF LEMMA 5

To simplify the notation, let $\mathbf{S} := \mathbf{CE}^{-1}\mathbf{C}^{\top}$. Using Lemma 2, we arrive at

$$\mathbf{Dx}^{k+1} = -\frac{1}{\rho}\nabla F(\mathbf{x}^{k+1}) + \mathbf{Sx}^{k} - (\mathbf{D}-\mathbf{S})\sum_{t=0}^{k}\mathbf{x}^{t}.$$

Subtracting $(\mathbf{D}-\mathbf{S})\mathbf{x}^{k+1}$ from both sides yields

$$\mathbf{S}(\mathbf{x}^{k+1}-\mathbf{x}^{k}) = -\frac{1}{\rho}\nabla F(\mathbf{x}^{k+1}) - (\mathbf{D}-\mathbf{S})\sum_{t=0}^{k+1}\mathbf{x}^{t}.$$

Noticing that $\mathbf{Q} = (\mathbf{D}-\mathbf{S})^{1/2}$ and $\mathbf{r}^{k} = \sum_{t=0}^{k}\mathbf{Q}\mathbf{x}^{t}$, we obtain

$$\mathbf{S}(\mathbf{x}^{k+1}-\mathbf{x}^{k}) = -\frac{1}{\rho}\nabla F(\mathbf{x}^{k+1}) - \mathbf{Qr}^{k+1}. \qquad (43)$$

The proof of Lemma 2 shows that if $\mathbf{y}^{0}=\mathbf{0}$, then $\mathbf{y}^{k}=\rho\mathbf{Qr}^{k}$, which leads to $\mathbf{y}^{\star}=\rho\mathbf{Qr}^{\star}$ as $k\rightarrow\infty$. Using (25), we arrive at

$$\nabla F(\mathbf{x}^{\star}) = -\rho\mathbf{Qr}^{\star}.$$

Combining this with (43) completes the proof.

## APPENDIX F
## PROOF OF THEOREM 1

It suffices to prove that

$$\|\mathbf{q}^{k+1}-\mathbf{q}^{\star}\|_{\mathbf{G}}^{2} \leq \frac{1}{1+\delta}\|\mathbf{q}^{k}-\mathbf{q}^{\star}\|_{\mathbf{G}}^{2}. \qquad (44)$$

Lemma 5 and the strong convexity of $F(\cdot)$ lead to

$$\frac{2}{\rho}\sigma\|\mathbf{x}^{k+1}-\mathbf{x}^{\star}\|_{2}^{2}$$
$$\leq \frac{2}{\rho}(\mathbf{x}^{k+1}-\mathbf{x}^{\star})^{\top}(\nabla F(\mathbf{x}^{k+1})-\nabla F(\mathbf{x}^{\star}))$$
$$= -2(\mathbf{r}^{k+1}-\mathbf{r}^{\star})^{\top}\mathbf{Q}(\mathbf{r}^{k+1}-\mathbf{r}^{\star})$$
$$\quad -2(\mathbf{x}^{k+1}-\mathbf{x}^{k})^{\top}\mathbf{CE}^{-1}\mathbf{C}^{\top}(\mathbf{x}^{k+1}-\mathbf{x}^{\star})$$
$$= 2(\mathbf{q}^{k+1}-\mathbf{q}^{k})^{\top}\mathbf{G}(\mathbf{r}^{k+1}-\mathbf{r}^{\star})$$
$$= \|\mathbf{q}^{k}-\mathbf{q}^{\star}\|_{\mathbf{G}}^{2}-\|\mathbf{q}^{k+1}-\mathbf{q}^{\star}\|_{\mathbf{G}}^{2}-\|\mathbf{q}^{k}-\mathbf{q}^{k+1}\|_{\mathbf{G}}^{2}.$$

Similarly, Lemma 5 and strong convexity imply that

$$\frac{2}{\rho}\frac{1}{L}\|\nabla F(\mathbf{x}^{k+1})-\nabla F(\mathbf{x}^{\star})\|_{2}^{2}$$
$$\leq \|\mathbf{q}^{k}-\mathbf{q}^{\star}\|_{\mathbf{G}}^{2}-\|\mathbf{q}^{k+1}-\mathbf{q}^{\star}\|_{\mathbf{G}}^{2}-\|\mathbf{q}^{k}-\mathbf{q}^{k+1}\|_{\mathbf{G}}^{2}.$$

For any $\beta\in(0,1)$, we have

$$\beta\frac{2}{\rho}\sigma\|\mathbf{x}^{k+1}-\mathbf{x}^{\star}\|_{2}^{2} + (1-\beta)\frac{2}{\rho}\frac{1}{L}\|\nabla F(\mathbf{x}^{k+1})-\nabla F(\mathbf{x}^{\star})\|_{2}^{2}$$
$$\leq \|\mathbf{q}^{k}-\mathbf{q}^{\star}\|_{\mathbf{G}}^{2}-\|\mathbf{q}^{k+1}-\mathbf{q}^{\star}\|_{\mathbf{G}}^{2}-\|\mathbf{q}^{k}-\mathbf{q}^{k+1}\|_{\mathbf{G}}^{2}.$$

To prove (44), it suffices to show that

$$\|\mathbf{q}^{k}-\mathbf{q}^{k+1}\|_{\mathbf{G}}^{2} + \beta\frac{2}{\rho}\sigma\|\mathbf{x}^{k+1}-\mathbf{x}^{\star}\|_{2}^{2}$$
$$+ (1-\beta)\frac{2}{\rho}\frac{1}{L}\|\nabla F(\mathbf{x}^{k+1}-\nabla F(\mathbf{x}^{\star}))\|_{2}^{2}$$
$$\geq \delta\|\mathbf{q}^{k+1}-\mathbf{q}^{\star}\|_{\mathbf{G}}^{2}$$

which is equivalent to

$$(1-\beta)\frac{2}{\rho}\frac{1}{L}\|\nabla F(\mathbf{x}^{k+1})-\nabla F(\mathbf{x}^{\star})\|_{2}^{2} + \|\mathbf{q}^{k}-\mathbf{q}^{k+1}\|_{\mathbf{G}}^{2}$$
$$+ \|\mathbf{x}^{k+1}-\mathbf{x}^{\star}\|_{\mathbf{M}}^{2} \geq \delta\|\mathbf{r}^{k+1}-\mathbf{r}^{\star}\|_{2}^{2} \quad (45)$$

where $\mathbf{M} := \frac{2\sigma\beta}{\rho}\mathbf{I}-\delta\mathbf{CE}^{-1}\mathbf{C}^{\top}$. Observing the left hand side of (45), it suffices to show that

$$\|\mathbf{x}^{k+1}-\mathbf{x}^{\star}\|_{\mathbf{M}}^{2} + (1-\beta)\frac{2}{\rho L}\|\nabla F(\mathbf{x}^{k+1})-\nabla F(\mathbf{x}^{\star})\|_{2}^{2}$$
$$\geq \delta\|\mathbf{r}^{k+1}-\mathbf{r}^{\star}\|_{2}^{2}. \quad (46)$$

Since $\mathbf{r}^{k+1}$ and $\mathbf{r}^{\star}$ are orthogonal to $\mathbf{1}$ and the null space of $\mathbf{Q}$ is span$\{\mathbf{1}\}$, both vectors belong to the column space of $\mathbf{Q}$. Using Lemma 4, we obtain

$$\delta\|\mathbf{r}^{k+1}-\mathbf{r}^{\star}\|_{2}^{2} \leq \frac{\delta}{\lambda}\|\mathbf{Q}(\mathbf{r}^{k+1}-\mathbf{r}^{\star})\|_{2}^{2}$$
$$\leq \frac{\delta}{\lambda}\|\mathbf{M}(\mathbf{x}^{k+1}-\mathbf{x}^{\star}) - \frac{1}{\rho}(\nabla F(\mathbf{x}^{k+1})-\nabla F(\mathbf{x}^{\star}))\|_{2}^{2}$$
$$\leq \frac{2\delta\Lambda}{\lambda}\|\mathbf{x}^{k+1}-\mathbf{x}^{\star}\|_{\mathbf{M}}^{2} + \frac{2\delta}{\rho\lambda}\|\nabla F(\mathbf{x}^{k+1})-\nabla F(\mathbf{x}^{\star})\|_{2}^{2}.$$
$$(47)$$

Comparing (47) with (46) suggests that it is sufficient to have

$$\delta \leq \min\left\{\frac{2\beta\sigma}{\rho(\Lambda+\frac{2\Lambda}{\lambda})}, \frac{(1-\beta)\rho\lambda}{L}\right\}. \qquad (48)$$

## REFERENCES

[1] J. A. Bazerque and G. B. Giannakis. Distributed Spectrum Sensing for Cognitive Radio Networks by Exploiting Sparsity. *IEEE Trans. Signal Process.*, 58(3):1847–1862, March 2010.

[2] Juan Andrés Bazerque, Gonzalo Mateos, and Georgios B. Giannakis. Group-lasso on splines for spectrum cartography. *Signal Process. IEEE Trans. On*, 59(10):4648–4663, 2011.

[3] Dimitri P. Bertsekas and John N. Tsitsiklis. *Parallel and Distributed Computation: Numerical Methods*, volume 23. Prentice-Hall, Inc., Upper Saddle River, NJ, USA, 1989.

[4] Daniel Boley. Local linear convergence of the alternating direction method of multipliers on quadratic or linear programs. *SIAM J. Optim.*, 23(4):2183–2207, 2013.

[5] Stephen Boyd, Neal Parikh, Eric Chu, Borja Peleato, and Jonathan Eckstein. Distributed optimization and statistical learning via the alternating direction method of multipliers. *Found. Trends® Mach. Learn.*, 3(1):1–122, 2011.

[6] Fan R. K. Chung. *Spectral Graph Theory.* Number no. 92. American Mathematical Society, Providence, R.I, 1997.

[7] Wei Deng and Wotao Yin. On the Global and Linear Convergence of the Generalized Alternating Direction Method of Multipliers. *J Sci Comput*, 66(3):889–916, May 2015.

[8] A. G. Dimakis, S. Kar, J. M. F. Moura, M. G. Rabbat, and A. Scaglione. Gossip Algorithms for Distributed Signal Processing. *Proc. IEEE*, 98(11):1847–1864, November 2010.

[9] J. C. Duchi, A. Agarwal, and M. J. Wainwright. Dual Averaging for Distributed Optimization: Convergence Analysis and Network Scaling. *IEEE Trans. Autom. Control*, 57(3):592–606, March 2012.

[10] Jonathan Eckstein and Dimitri P. Bertsekas. On the Douglas—Rachford splitting method and the proximal point algorithm for maximal monotone operators. *Mathematical Programming*, 55(1-3):293–318, April 1992.

[11] T. Erseghe, D. Zennaro, E. Dall'Anese, and L. Vangelista. Fast Consensus by the Alternating Direction Multipliers Method. *IEEE Trans. Signal Process.*, 59(11):5523–5537, November 2011.

[12] Michalis Faloutsos, Petros Faloutsos, and Christos Faloutsos. On power-law relationships of the internet topology. In *ACM SIGCOMM computer communication review*, volume 29, pages 251–262. ACM, 1999.

[13] P. A. Forero, A. Cano, and G. B. Giannakis. Distributed Clustering Using Wireless Sensor Networks. *IEEE J. Sel. Top. Signal Process.*, 5(4):707–724, August 2011.

[14] Pedro A. Forero, Alfonso Cano, and Georgios B. Giannakis. Consensus-Based Distributed Support Vector Machines. *J. Mach. Learn. Res.*, 11(May):1663–1707, 2010.

[15] D. Gabay. Applications of the Method of Multipliers to Variational Inequalities. In Michel Fortin and Roland Glowinski, editors, *Studies in Mathematics and Its Applications*, volume 15 of *Augmented Lagrangian Methods: Applications to the Numerical Solution of Boundary-Value Problems*, pages 299–331. North-Holland, Amsterdam, January 1983.

[16] Georgios B. Giannakis, Qing Ling, Gonzalo Mateos, Ioannis D. Schizas, and Hao Zhu. Decentralized Learning for Wireless Communications and Networking. In Roland Glowinski, Stanley J. Osher, and Wotao Yin, editors, *Splitting Methods in Communication, Imaging, Science, and Engineering*, pages 461–497. Springer International Publishing, Cham, 2016.

[17] Roland Glowinski and A. Marroco. Sur l'approximation, par éléments finis d'ordre un, et la résolution, par pénalisation-dualité d'une classe de problèmes de Dirichlet non linéaires. *Rev. Fr. Autom. Inform. Rech. Opérationnelle Anal. Numér.*, 9(R2):41–76, 1975.

[18] Tom Goldstein, Brendan O'Donoghue, Simon Setzer, and Richard Baraniuk. Fast alternating direction optimization methods. *SIAM J. Imaging Sci.*, 7(3):1588–1623, 2014.

[19] Jason P Gross, Ranjit B Pandit, Clive G Cook, and Thomas H Matson. Network with distributed shared memory, April 19 2016. US Patent 9,317,469.

[20] Mingyi Hong and Zhi-Quan Luo. On the linear convergence of the alternating direction method of multipliers. *Math. Program.*, 162(1-2):165–199, March 2017.

[21] F. Iutzeler, P. Bianchi, P. Ciblat, and W. Hachem. Explicit Convergence Rate of a Distributed Alternating Direction Method of Multipliers. *IEEE Trans. Autom. Control*, 61(4):892–904, April 2016.

[22] Anil K. Jain. *Fundamentals of Digital Image Processing*. Prentice-Hall, Inc., Upper Saddle River, NJ, USA, 1989.

[23] Vassilis Kekatos and Georgios B. Giannakis. Distributed robust power system state estimation. *IEEE Trans. Power Syst.*, 28(2):1617–1626, 2013.

[24] Q. Ling and Z. Tian. Decentralized Sparse Signal Recovery for Compressive Sleeping Wireless Sensor Networks. *IEEE Trans. Signal Process.*, 58(7):3816–3827, July 2010.

[25] Qing Ling, Yaohua Liu, Wei Shi, and Zhi Tian. Communication-efficient weighted ADMM for decentralized network optimization. In *2016 IEEE International Conference on Acoustics, Speech and Signal Processing*, pages 4821–4825, 2016.

[26] Meng Ma, Athanasios N. Nikolakopoulos, and Georgios B. Giannakis. Fast Decentralized Learning via Hybrid Consensus ADMM. In *2018 IEEE International Conference on Acoustics, Speech and Signal Processing*, April 2018.

[27] A. Makhdoumi and A. Ozdaglar. Convergence Rate of Distributed ADMM Over Networks. *IEEE Trans. Autom. Control*, 62(10):5082–5095, October 2017.

[28] G. Mateos, J. A. Bazerque, and G. B. Giannakis. Distributed Sparse Linear Regression. *IEEE Trans. Signal Process.*, 58(10):5262–5276, October 2010.

[29] A. Nedic and A. Ozdaglar. Distributed Subgradient Methods for Multi-Agent Optimization. *IEEE Trans. Autom. Control*, 54(1):48–61, January 2009.

[30] Angelia Nedic, Asuman Ozdaglar, and Pablo A. Parrilo. Constrained consensus and optimization in multi-agent networks. *IEEE Trans. Autom. Control*, 55(4):922–938, 2010.

[31] Mark EJ Newman. The structure and function of complex networks. *SIAM review*, 45(2):167–256, 2003.

[32] Athanasios N. Nikolakopoulos and John D. Garofalakis. On the multi-level near complete decomposability of a class of multiprocessing systems. In *Proceedings of the 2016 PCI Conference*, PCI '16, pages 14:1–14:6, New York, NY, USA, 2016. ACM.

[33] R. Olfati-Saber. Distributed Kalman Filter with Embedded Consensus Filters. In *Proceedings of the 44th IEEE Conference on Decision and Control*, pages 8179–8184, December 2005.

[34] R. Olfati-Saber and J. S. Shamma. Consensus Filters for Sensor Networks and Distributed Sensor Fusion. In *Proceedings of the 44th IEEE Conference on Decision and Control*, pages 6698–6703, December 2005.

[35] J. B. Predd, S. R. Kulkarni, and H. V. Poor. A Collaborative Training Algorithm for Distributed Learning. *IEEE Trans. Inf. Theory*, 55(4):1856–1871, April 2009.

[36] Wei Ren. Consensus based formation control strategies for multi-vehicle systems. In *2006 American Control Conference*, pages 4237–4242, June 2006.

[37] Wei Ren, Randal W. Beard, and Ella M. Atkins. Information consensus in multivehicle cooperative control. *IEEE Control Syst. Mag.*, 27(2):71–82, April 2007.

[38] I. D. Schizas, A. Ribeiro, and G. B. Giannakis. Consensus in Ad Hoc WSNs With Noisy Links – Part I: Distributed Estimation of Deterministic Signals. *IEEE Trans. Signal Process.*, 56(1):350–364, January 2008.

[39] Wei Shi, Qing Ling, Kun Yuan, Gang Wu, and Wotao Yin. On the Linear Convergence of the ADMM in Decentralized Consensus Optimization. *IEEE Trans. Signal Process.*, 62(7):1750–1761, April 2014.

[40] S. Sundhar Ram, A. Nedić, and V. V. Veeravalli. Distributed Stochastic Subgradient Projection Algorithms for Convex Optimization. *J Optim Theory Appl*, 147(3):516–545, December 2010.

[41] Hakan Terelius, Ufuk Topcu, and Richard M. Murray. Decentralized multi-agent optimization via dual decomposition. *IFAC Proc. Vol.*, 44(1):11245–11251, 2011.

[42] K. I. Tsianos, S. Lawlor, and M. G. Rabbat. Push-Sum Distributed Dual Averaging for convex optimization. In *2012 IEEE 51st IEEE Conference on Decision and Control (CDC)*, pages 5453–5458, December 2012.

[43] K. I. Tsianos and M. G. Rabbat. Distributed strongly convex optimization. In *2012 50th Annual Allerton Conference on Communication, Control, and Computing (Allerton)*, pages 593–600, October 2012.

[44] J. N. Tsitsiklis. *Problems in Decentralized Decision Making and Computation*. PhD thesis, Massachusetts Institute of Technology, December 1984.

[45] Lin Xiao, Stephen Boyd, and Sanjay Lall. A Scheme for Robust Distributed Sensor Fusion Based on Average Consensus. In *Proceedings of the 4th International Symposium on Information Processing in Sensor Networks*, IPSN '05, pages 63–70, April 2005.

[46] Hao Zhu, Alfonso Cano, and Georgios Giannakis. Distributed consensus-based demodulation: Algorithms and error analysis. *IEEE Trans. Wirel. Commun.*, 9(6):2044–2054, June 2010.